\input amstex
\documentstyle{amsppt}
%
\catcode`@=11
\redefine\output@{%
  \def\break{\penalty-\@M}\let\par\endgraf
  \ifodd\pageno\global\hoffset=105pt\else\global\hoffset=8pt\fi  
  \shipout\vbox{%
    \ifplain@
      \let\makeheadline\relax \let\makefootline\relax
    \else
      \iffirstpage@ \global\firstpage@false
        \let\rightheadline\frheadline
        \let\leftheadline\flheadline
      \else
        \ifrunheads@ 
        \else \let\makeheadline\relax
        \fi
      \fi
    \fi
    \makeheadline \pagebody \makefootline}%
  \advancepageno \ifnum\outputpenalty>-\@MM\else\dosupereject\fi
}
\def\Beta{\mathchar"0\hexnumber@\rmfam 42}
\catcode`\@=\active
\nopagenumbers
\def\negskp{\hskip -2pt}

\def\blue#1{#1}

\catcode`#=11\def\diez{#}\catcode`#=6
\catcode`&=11\catcode`&=4
\catcode`_=11\catcode`_=8
\def\mycite#1{\cite{\blue{#1}}\immediate\special{ps:
     ShrHPSdict begin /ShrBORDERthickness 0 def}}
\def\myciterange#1#2#3#4{\cite{\blue{#2#3#4}}\immediate\special{ps:
     ShrHPSdict begin /ShrBORDERthickness 0 def}}
\def\mytag#1{%
    \tag#1}
\def\mythetag#1{\thetag{\blue{#1}}\immediate\special{ps:
     ShrHPSdict begin /ShrBORDERthickness 0 def}}
\def\myrefno#1{\no#1}
\def\myhref#1#2{\blue{#2}\immediate\special{ps:
     ShrHPSdict begin /ShrBORDERthickness 0 def}}
\def\myEarXivlink{\myhref{http://arXiv.org}{http:/\negskp/arXiv.org}}

\def\mytheorem#1{\csname proclaim\endcsname{Theorem #1}}
\def\mytheoremwithtitle#1#2{\csname proclaim\endcsname{Theorem #1#2}}
\def\mythetheorem#1{\blue{#1}\immediate\special{ps:
     ShrHPSdict begin /ShrBORDERthickness 0 def}}
\def\mylemma#1{\csname proclaim\endcsname{Lemma #1}}
\def\mylemmawithtitle#1#2{\csname proclaim\endcsname{Lemma #1#2}}

\def\mycorollary#1{\csname proclaim\endcsname{Corollary #1}}

\def\myconjecture#1{\csname proclaim\endcsname{Conjecture #1}}
\def\myconjecturewithtitle#1#2{\csname proclaim\endcsname{Conjecture #1#2}}
\def\mytheconjecture#1{\blue{#1}\immediate\special{ps:
     ShrHPSdict begin /ShrBORDERthickness 0 def}}

\pagewidth{360pt}
\pageheight{606pt}
\topmatter
\title
Perfect cuboids and irreducible polynomials.
\endtitle
\author
Ruslan Sharipov
\endauthor
\address Bashkir State University, 32 Zaki Validi street, 450074 Ufa, Russia
\endaddress
\email\myhref{mailto:r-sharipov\@mail.ru}{r-sharipov\@mail.ru}
\endemail
\abstract
    The problem of constructing a perfect cuboid is related to a certain class 
of univariate polynomials with three integer parameters $a$, $b$, and $u$. Their 
irreducibility over the ring of integers under certain restrictions for $a$, $b$, 
and $u$ would mean the non-existence of perfect cuboids.  This irreducibility 
is conjectured and then verified numerically for approximately $10\,000$ instances 
of $a$, $b$, and $u$. 
\endabstract
\subjclassyear{2000}
\subjclass 11D41, 11D72, 11Y50, 12E05\endsubjclass
\endtopmatter
\TagsOnRight
\document

\head
1. Introduction.
\endhead
     An Euler cuboid is a rectangular parallelepiped whose edges and 
face diagonals all have integer lengths. A perfect cuboid is an Euler 
cuboid whose space diagonal is also of an integer length. Cuboids with integer 
sides and face diagonals were known before Euler (see \mycite{1} and \mycite{2}). 
However, they became famous due to Leonhard Euler (see \mycite{3}) and were 
named after him.\par
     As for perfect cuboids, none of them is known by now. The problem of finding
perfect cuboids or proving their non-existence is an open mathematical problem.
The search for perfect cuboids has the long history. It is reflected in 
\myciterange{4}{4}{--}{34}.\par
     In \mycite{35} the problem of finding a perfect cuboid was reduced to the 
following Diophantine equation of the order 12 with four variables $a$, $b$, $c$,
and $u$:
$$
\gathered
u^4\,a^4\,b^4+6\,a^4\,u^2\,b^4\,c^2-2\,u^4\,a^4\,b^2\,c^2-2\,u^4\,a^2
\,b^4\,c^2+4\,u^2\,b^4\,a^2\,c^4+\\
+\,4\,a^4\,u^2\,b^2\,c^4-12\,u^4\,a^2\,b^2\,c^4+u^4\,a^4\,c^4+u^4\,b^4\,c^4
+a^4\,b^4\,c^4+\\
+\,6\,a^4\,u^2\,c^6+6\,u^2\,b^4\,c^6-8\,a^2\,b^2\,u^2\,c^6-2\,u^4\,a^2\,c^6
-2\,u^4\,b^2\,c^6-\\
-\,2\,a^4\,b^2\,c^6-2\,b^4\,a^2\,c^6+u^4\,c^8+b^4\,c^8+a^4\,c^8
+4\,a^2\,u^2\,c^8+\\
+\,4\,b^2\,u^2\,c^8-12\,b^2\,a^2\,c^8+6\,u^2\,c^{10}-2\,a^2\,c^{10}
-2\,b^2\,c^{10}+c^{12}=0.
\endgathered
\quad
\mytag{1.1}
$$
The exact result of the paper \mycite{35} is formulated as follows.
\mytheorem{1.1} A perfect Euler cuboid does exist if and only if 
the Diophantine equation \mythetag{1.1} has a solution such that $a$, $b$,
$c$, and $u$ are positive integer numbers obeying the inequalities 
and $a<c$, $b<c$, $u<c$, and $(a+c)\,(b+c)>2\,c^2$.
\endproclaim
     A more simple equation associated with perfect cuboids was derived in
\mycite{18} (see also \mycite{27}). However, our goal in this paper is to 
study the equation \mythetag{1.1} (because it is new) and derive the results 
declared in the abstract. 
\head
2. Rational cuboids.
\endhead
     A rational cuboid is a rectangular parallelepiped the lengths of whose 
edges are rational numbers. If the lengths of face diagonals are also rational
numbers, it is called a rational Euler cuboid. Finally, if the 
length of the space diagonal is a rational number too, we have a perfect rational 
cuboid. It is easy to see that each rational Euler cuboid can be transformed to 
an Euler cuboid with integer sides and diagonals. In the case of perfect 
cuboids (either integer or rational) each such cuboid can be transformed to a 
perfect rational cuboid whose space diagonal is equal to unity (see \mycite{35}). 
Conversely each perfect rational cuboid with unit space diagonal yields some 
perfect cuboid with integer sides and diagonals. Therefore, saying a perfect 
rational cuboid, we assume its space diagonal to be equal to unity. 
\head
3. Expressions for the sides and face diagonals.
\endhead
     Note that the equation \mythetag{1.1} is homogeneous with respect to its
variables $a$, $b$, $c$, and $u$. Since $c>0$ in the theorem~\mythetheorem{1.1},
we can introduce the fractions
$$
\xalignat 3
&\hskip -2em
\alpha=\frac{a}{c},
&&\beta=\frac{b}{c},
&&\upsilon=\frac{u}{c}.
\mytag{3.1}
\endxalignat
$$
In terms of the rational variables \mythetag{3.1} the equation \mythetag{1.1}
is written as 
$$
\gathered
\upsilon^4\,\alpha^4\,\beta^4+(6\,\alpha^4\,\upsilon^2\,\beta^4
-2\,\upsilon^4\,\alpha^4\,\beta^2-2\,\upsilon^4\,\alpha^2
\,\beta^4)+(4\,\upsilon^2\,\beta^4\,\alpha^2+\\
+\,4\,\alpha^4\,\upsilon^2\,\beta^2-12\,\upsilon^4\,\alpha^2\,\beta^2
+\upsilon^4\,\alpha^4+\upsilon^4\,\beta^4
+\alpha^4\,\beta^4)+(6\,\alpha^4\,\upsilon^2+6\,\upsilon^2\,\beta^4-\\
-\,8\,\alpha^2\,\beta^2\,\upsilon^2-2\,\upsilon^4\,\alpha^2
-2\,\upsilon^4\,\beta^2-2\,\alpha^4\,\beta^2-2\,\beta^4\,\alpha^2)
+(\upsilon^4+\beta^4+\\
+\,\alpha^4+4\,\alpha^2\,\upsilon^2+4\,\beta^2\,\upsilon^2
-12\,\beta^2\,\alpha^2)+(6\,\upsilon^2-2\,\alpha^2-2\,\beta^2)+1=0.
\endgathered
\quad
\mytag{3.2}
$$\par
     Note that the variables $a$, $b$, $c$, and $u$ in \mythetag{1.1} are
neither edges nor diagonals of a perfect cuboid, they are just parameters. 
They yield the rational parameters $\alpha$, $\beta$, and $\upsilon$ in
\mythetag{3.2} according to the formulas \mythetag{3.1}. The edges and 
face diagonals of a perfect rational cuboid are expressed through $\alpha$, 
$\beta$, and $\upsilon$. Let's denote through $x_1$, $x_2$, and $x_3$ the 
edges of such a cuboid and through $d_1$, $d_2$, and $d_3$ its side diagonals:
$$
\kern -10pt
\!\!(x_1)^2+(x_2)^2=(d_3)^2\!,\ \ \quad
(x_2)^2+(x_3)^2=(d_1)^2\!,\ \ \quad
(x_3)^2+(x_1)^2=(d_2)^2\!.
\mytag{3.3}
$$
Then $x_1$ and $d_1$ are expressed through the parameter $\upsilon$:
$$
\xalignat 2
&\hskip -2em 
x_1=\frac{2\,\upsilon}{1+\upsilon^2},
&&d_1=\frac{1-\upsilon^2}{1+\upsilon^2}.
\mytag{3.4}
\endxalignat
$$
Let's denote through $z$ the following auxiliary parameter:
$$
\hskip -2em
z=\frac{(1+\upsilon^2)\,(1-\beta^2)\,(1+\alpha^2)}{2\,(1+\beta^2)\,(1
-\alpha^2\,\upsilon^2)}.
\mytag{3.5}
$$
Then the edges $x_2$ and $x_3$ are expressed by the formulas 
$$
\xalignat 2
&\hskip -2em
x_2=\frac{2\,z\,(1-\upsilon^2)}{(1+\upsilon^2)\,(1+z^2)},
&&x_3=\frac{(1-\upsilon^2)\,(1-z^2)}{(1+\upsilon^2)\,(1+z^2)}.
\mytag{3.6}
\endxalignat
$$
The side diagonals $d_2$ and $d_3$ are given by the following formulas:
$$
\hskip -2em
\gathered
d_2=\frac{(1+\upsilon^2)\,(1+z^2)+2\,z(1-\upsilon^2)}{(1+\upsilon^2)\,(1+z^2)}
\,\beta,\\
\vspace{1ex}
d_3=\frac{2\,(\upsilon^2\,z^2+1)}{(1+\upsilon^2)\,(1+z^2)}\,\alpha.
\endgathered
\mytag{3.7}
$$\par
     The formulas \mythetag{3.4}, \mythetag{3.5}, \mythetag{3.6} and 
\mythetag{3.7} are taken from \mycite{35}. They can be verified by means of
direct calculations. Indeed, the second equality \mythetag{3.3} turns to
an identity due to the formulas \mythetag{3.6}. Apart from the equalities
\mythetag{3.3}, a perfect rational cuboid is characterized by the equalities
$$
\xalignat 3
&\hskip -2em
(x_1)^2+(d_1)^2=1,
&&(x_2)^2+(d_2)^2=1,
&&(x_3)^2+(d_3)^2=1.
\qquad
\mytag{3.8}
\endxalignat
$$
They mean that the space diagonal of such a cuboid is equal to unity. The 
first equality \mythetag{3.8} turns to an identity due to the formulas 
\mythetag{3.4}.\par
    Thus, the second equality \mythetag{3.3} and the first equality 
\mythetag{3.8} turn to identities. Other four equalities in \mythetag{3.3}
and \mythetag{3.8} also turn to identities due to \mythetag{3.4}, 
\mythetag{3.5}, \mythetag{3.6} and \mythetag{3.7} modulo the equation
\mythetag{3.2}. 
\head
4. Back to integer numbers. 
\endhead
     The equation \mythetag{1.1} is homogeneous with respect to its variables. 
Therefore due to \mythetag{3.1} and due to the theorem~\mythetheorem{1.1} the 
parameters $a$, $b$, $c$, and $u$ in the equation \mythetag{1.1} can be treated 
as a quadruple of positive coprime integer numbers, i\.\,e\. their greatest 
common divisor is equal to unity:
$$
\hskip -2em
\gcd(a,b,c,u)=1.
\mytag{4.1}
$$
Let's denote through $m$ the greatest common divisor of $a$, $b$, and $u$:
$$
\hskip -2em
\gcd(a,b,u)=m.
\mytag{4.2}
$$
Then from \mythetag{4.1} and \mythetag{4.2} we derive the equality 
$$
\hskip -2em
\gcd(m,c)=1,
\mytag{4.3}
$$
i\.\,e\. $m$ and $c$ are coprime. Due to \mythetag{4.2} and \mythetag{4.3} 
the fractions $a/m$, $b/m$, and $u/m$ reduce to integer numbers, while $c/m$ 
is an irreducible fraction if $m\neq 1$. The formula \mythetag{3.1} can be 
written in terms of these fractions:
$$
\xalignat 3
&\hskip -2em
\alpha=\frac{a/m}{c/m},
&&\beta=\frac{b/m}{c/m},
&&\upsilon=\frac{u/m}{c/m}.
\mytag{4.4}
\endxalignat
$$
Relying on \mythetag{4.4}, we can change variables as follows:
$$
\xalignat 4
&\hskip -2em
\frac{a}{m}\to a,
&&\frac{b}{m}\to b,
&&\frac{u}{m}\to u,
&&\frac{c}{m}\to t.
\quad
\mytag{4.5}
\endxalignat 
$$
In terms of the new variable $t=c/m$ and in terms of the renewed variables $a$, $b$, 
and $u$ in \mythetag{4.5} the formulas \mythetag{4.4} are written as 
$$
\xalignat 3
&\hskip -2em
\alpha=\frac{a}{t},
&&\beta=\frac{b}{t},
&&\upsilon=\frac{u}{t},
\mytag{4.6}
\endxalignat
$$
while the equation \mythetag{1.1} turns to the following equation:
$$
\gathered
t^{12}+(6\,u^2\,-2\,a^2\,-2\,b^2)\,t^{10}
+(u^4\,+b^4+a^4+4\,a^2\,u^2+\\
+\,4\,b^2\,u^2-12\,b^2\,a^2)\,t^8+(6\,a^4\,u^2+6\,u^2\,b^4-8\,a^2\,b^2\,u^2-\\
-\,2\,u^4\,a^2-2\,u^4\,b^2-2\,a^4\,b^2-2\,b^4\,a^2)\,t^6+(4\,u^2\,b^4\,a^2+\\
+\,4\,a^4\,u^2\,b^2-12\,u^4\,a^2\,b^2+u^4\,a^4+u^4\,b^4+a^4\,b^4)\,t^4+\\
+\,(6\,a^4\,u^2\,b^4-2\,u^4\,a^4\,b^2-2\,u^4\,a^2
\,b^4)\,t^2+u^4\,a^4\,b^4=0.
\endgathered
\quad
\mytag{4.7}
$$
As for the formula \mythetag{4.2}, for the renewed variables $a$, $b$, and $u$ in
\mythetag{4.5} it yields
$$
\hskip -2em
\gcd(a,b,u)=1.
\mytag{4.8}
$$
The formula \mythetag{4.8} means that $a$, $b$, and $u$ in \mythetag{4.6} and 
\mythetag{4.7} are coprime.\par
     Note that the equation \mythetag{4.7} is the same as the initial equation
\mythetag{1.1}, but $c$ is replaced by $t$ and the terms are reordered like in 
a univariate polynomial of the variable $t$. The theorem~\mythetheorem{1.1} now
is reformulated as follows.
\mytheorem{4.1} A perfect Euler cuboid does exist if and only if for some
positive coprime integer numbers $a$, $b$, and $u$ the polynomial equation 
\mythetag{4.7} has a rational solution $t$ obeying the inequalities
$t>a$, $t>b$, $t>u$, and $(a+t)\,(b+t)>2\,t^2$. 
\endproclaim
\head
5. Factoring the polynomial equation. 
\endhead
     Let's denote through $P_{abu}(t)$ the polynomial in the left hand side
of the equation \mythetag{4.7}. Denoting it in this way, we shall treat it
as a univariate polynomial of $t$, while the variables $a$, $b$, and $u$ are 
treated as parameters:
$$
\gathered
P_{abu}(t)=t^{12}+(6\,u^2\,-2\,a^2\,-2\,b^2)\,t^{10}
+(u^4\,+b^4+a^4+4\,a^2\,u^2+\\
+\,4\,b^2\,u^2-12\,b^2\,a^2)\,t^8+(6\,a^4\,u^2+6\,u^2\,b^4-8\,a^2\,b^2\,u^2-\\
-\,2\,u^4\,a^2-2\,u^4\,b^2-2\,a^4\,b^2-2\,b^4\,a^2)\,t^6+(4\,u^2\,b^4\,a^2+\\
+\,4\,a^4\,u^2\,b^2-12\,u^4\,a^2\,b^2+u^4\,a^4+u^4\,b^4+a^4\,b^4)\,t^4+\\
+\,(6\,a^4\,u^2\,b^4-2\,u^4\,a^4\,b^2-2\,u^4\,a^2
\,b^4)\,t^2+u^4\,a^4\,b^4.
\endgathered
\quad
\mytag{5.1}
$$
The polynomial \mythetag{5.1} is symmetric with respect to the parameters 
$a$ and $b$, i\.\,e\.
$$
\hskip -2em
P_{abu}(t)=P_{\kern 0.5pt bau}(t).
\mytag{5.2}
$$
In order to study the polynomial $P_{abu}(t)$ we consider some special cases:
$$
\pagebreak
\xalignat 3
&\hskip -2em
\text{1) \ }a=b; &&\text{3) \ }b\,u=a^2; &&\text{5) \ }a=u;
\qquad\\
\vspace{-1.5ex}
\mytag{5.3}\\
\vspace{-1.5ex}
&\hskip -2em
\text{2) \ }a=b=u; &&\text{4) \ }a\,u=b^2; &&\text{6) \ }b=u.
\qquad
\endxalignat
$$\par
     {\bf The special case $a=b$}. In this special case the polynomial
$P_{abu}(t)=P_{aau}(t)$ is given by the following formula:
$$
\hskip -2em
\gathered
P_{aau}(t)=t^{12}+(6\,u^2-4\,a^2)\,t^{10}+(8\,a^2\,u^2-10\,a^4+u^4)\,t^8+\\
+\,(4\,a^4\,u^2-4\,a^6-4\,u^4\,a^2)\,t^6+(8\,a^6\,u^2+a^8-10\,u^4\,a^4)\,t^4+\\
+\,(6\,a^8\,u^2-4\,u^4\,a^6)\,t^2+u^4\,a^8.
\endgathered
\mytag{5.4}
$$
The polynomial \mythetag{5.4} is reducible. It is factored as
$$
\hskip -2em
P_{aau}(t)=(t^2+a^2)^2\,P_{au}(t),
\mytag{5.5}
$$
where the polynomial $P_{au}(t)$ is given by the formula
$$
\hskip -2em
\gathered
P_{au}(t)=t^8+6\,(u^2-a^2)\,t^6+(a^4-4\,a^2\,u^2+u^4)\,t^4-\\
-\,6\,a^2\,u^2\,(u^2-a^2)\,t^2+u^4\,a^4.
\endgathered
\mytag{5.6}
$$
The formulas \mythetag{5.5} and \mythetag{5.6} are easily proved by
direct calculations.\par
    {\bf The special case $a=b=u$}. This case corresponds to $a=u$ in
\mythetag{5.6}. If $a=u$, the polynomial $P_{au}(t)=P_{aa}(t)$ is reducible:
$$
\hskip -2em
P_{aa}(t)=(t-a)^2\,(t+a)^2\,(t^2+a^2)^2.
\mytag{5.7}
$$
Due to the coprimality \mythetag{4.8} the special case $a=b=u$ 
can fit the theorem~\mythetheorem{4.1} only for $a=b=u=1$. Then, due to 
\mythetag{5.5} and \mythetag{5.7}, the equation \mythetag{4.7} looks like
$$
\hskip -2em
(t-1)^2\,(t+1)^2\,(t^2+1)^4=0.
\mytag{5.8}
$$
The equation \mythetag{5.8} has two real rational solutions $t=-1$ and $t=1$.
Both of them do not fit the theorem~\mythetheorem{4.1}. Indeed, both of them do 
not satisfy the inequality $t>a$, where $a=1$.
\par
     Thus, the subcase $a=b=u$ of the special case $a=b$ do not provide any
perfect cuboid. Other subcases of the case $a=b$ are described by the following 
conjecture. 
\myconjecture{5.1} For any positive coprime integers $a\neq u$ the polynomial
$P_{au}(t)$ in \mythetag{5.6} is irreducible in the ring $\Bbb Z[t]$. 
\endproclaim
     {\bf The special case $b\,u=a^2$}. Combining $b\,u=a^2$ with
\mythetag{4.8} one can easily derive the following presentation
for the integer numbers $a$, $b$, and $u$:
$$
\xalignat 3
&\hskip -2em
a=p\,q,
&&b=p^2, &&u=q^2.
\mytag{5.9}
\endxalignat
$$
Here $p$ and $q$ are two positive integers, satisfying the coprimality
condition
$$
\pagebreak
\hskip -2em
\gcd(p,q)=1.
\mytag{5.10}
$$
Substituting \mythetag{5.9} into \mythetag{5.1}, we get the polynomial
$$
\gathered
P_{pqp^2\!q^2}(t)=t^{12}+(6\,q^4-2\,p^2\,q^2-2\,p^4)
\,t^{10}+(q^8+4\,p^2\,q^6+\\
+\,5\,p^4\,q^4-12\,p^6\,q^2+p^8)\,t^8
-2\,p^2\,q^2\,(q^8-2\,p^2\,q^6+4\,p^4\,q^4-\\
-\,2\,p^6\,q^2+p^8)\,t^6
+p^4\,q^4\,(q^8-12\,p^2\,q^6+5\,p^4\,q^4+4\,p^6\,q^2+p^8)\,t^4+\\
+\,q^8\,p^8\,(-2\,q^4-2\,p^2\,q^2+6\,p^4)\,t^2
+q^{12}\,p^{12}.
\endgathered
\quad
\mytag{5.11}
$$
The polynomial $P_{pqp^2\!q^2}(t)$ in \mythetag{5.11} is reducible. Indeed, 
we have
$$
\hskip -2em
P_{pqp^2\!q^2}(t)=(t-a)\,(t+a)\,Q_{pq}(t),
\mytag{5.12}
$$
where $Q_{pq}(t)$ is the following polynomial:
$$
\gathered
Q_{pq}(t)=t^{10}+(2\,q^2+p^2)\,(3\,q^2-2\,p^2)\,t^8
+(q^8+10\,p^2\,q^6+\\
+\,4\,p^4\,q^4-14\,p^6\,q^2+p^8)\,t^6
-p^2\,q^2\,(q^8-14\,p^2\,q^6+4\,p^4\,q^4+\\
+\,10\,p^6\,q^2+p^8)\,t^4-p^6\,q^6\,(q^2+2\,p^2)
\,(-2\,q^2+3\,p^2)\,t^2-q^{10}\,p^{10}.
\endgathered
\quad
\mytag{5.13}
$$
Due to \mythetag{5.12} the polynomial \mythetag{5.11} has two 
rational roots $t=a$ and $t=-a$. Both of them do not fit the 
theorem~\mythetheorem{4.1} since they do not satisfy the  
inequality $t>a$.\par 
    Other roots of the polynomial \mythetag{5.11}
coincide with the roots of the polynomial $Q_{pq}(t)$ in 
\mythetag{5.13}. The polynomial \mythetag{5.13} is reducible
if $q=p$. In this case we have
$$
\hskip -2em
Q_{pp}(t)=(t-a)\,(t+a)\,(t^2+a^2)^4.
\mytag{5.14}
$$
The formula \mythetag{5.14} is not surprising. For $q=p$ from
\mythetag{5.9} we derive $a=b=u$. This case was already considered
(see \mythetag{5.7} and \mythetag{5.8}). From $q=p$ and 
\mythetag{5.10} we derive $p=q=1$ and $a=b=u=1$.\par 
     In the case $p\neq q$ the polynomial \mythetag{5.13} is 
described by the following conjecture. 
\myconjecture{5.2} For any positive coprime integers $p\neq q$ the 
polynomial $Q_{pq}(t)$ in \mythetag{5.13} is irreducible in the ring 
$\Bbb Z[t]$. 
\endproclaim
     {\bf The special case $a\,u=b^2$}. This special case reduces to
the previous one. Indeed, from $a\,u=b^2$ and \mythetag{4.8} we derive
$$
\xalignat 3
&\hskip -2em
a=p^2,
&&b=p\,q, &&u=q^2,
\mytag{5.15}
\endxalignat
$$
where $p$ and $q$ are two positive integer numbers obeying the coprimality
condition \mythetag{5.10}. When substituted into \mythetag{5.1}, the formulas 
\mythetag{5.15} are equivalent to \mythetag{5.9} due to the symmetry
\mythetag{5.2}. They lead to the polynomial $P_{p^2\!pqq^2}(t)$ coinciding
with the polynomial \mythetag{5.11} and then lead to the polynomial 
\mythetag{5.13}, which was already considered.\par
     {\bf The special case $a=u$}. This special case is rather trivial. In this 
case the polynomial $P_{abu}(t)=P_{ubu}(t)$ in \mythetag{5.1} is reducible and 
we have the formula
$$
\hskip -2em
P_{ubu}(t)=(t^2+u^2)^4\,(t-b)^2\,(t+b)^2.
\mytag{5.16}
$$
The polynomial \mythetag{5.16} has two real rational roots $t=b$ and $t=-b$.
Both of them do not fit the theorem~\mythetheorem{4.1} \pagebreak since they 
do not satisfy the inequality $t>b$.\par
     {\bf The special case $b=u$}. This case is equivalent to the previous one
due to the symmetry \mythetag{5.2}.\par
     The general case not covered by the special cases listed in \mythetag{5.3} 
is described by the following conjecture.
\myconjecture{5.3} For any three positive coprime integer numbers $a$, $b$, and
$u$ such that none of the conditions \mythetag{5.3} is satisfied the polynomial 
\mythetag{5.1} is irreducible in the ring $\Bbb Z[t]$. 
\endproclaim
\head
6. Numeric study of the conjectures.
\endhead
     There are no proofs for the conjectures~\mytheconjecture{5.1}, 
\mytheconjecture{5.2}, and \mytheconjecture{5.3} at present time. 
Therefore I explored them numerically. For this purpose I used the Maxima 
package version 5.21.1 with the graphic interface wxMaxima 0.85 on the 
platform of Ubuntu 10.10 with Linux 2.5.35-24.\par
    The conjecture~\mytheconjecture{5.1} was verified and confirmed for $1\leqslant 
a\leqslant 100$ and $1\leqslant u\leqslant 100$. The conjecture~\mytheconjecture{5.2} 
was confirmed for $1\leqslant p\leqslant 100$ and $1\leqslant q\leqslant 100$. And 
the third conjecture~\mytheconjecture{5.3} was confirmed for $1\leqslant a\leqslant 22$,
$1\leqslant b\leqslant 22$, and $1\leqslant u\leqslant 22$. The number $22$ was
chosen intentionally since 
$$
22^3=10\,648\approx 10\,000=100^2.
$$
This equality means that each conjecture was tested and confirmed for approximately 
10\,000 instances of the numeric parameters in it. The overall result obtained can 
be formulated as follows: the equation \mythetag{4.7} has no solutions providing 
perfect cuboids for $a$, $b$, and $u$ less than or equal to $22$.
\head
7. Conclusions.
\endhead
     The conjectures~\mytheconjecture{5.1}, \mytheconjecture{5.2}, and 
\mytheconjecture{5.3} are not equivalent to the non-existence of perfect cuboids. 
However, if they are valid, this would be sufficient to prove that perfect cuboids
do not exist. The results of the numeric computations reported in the previous 
section support these conjectures. 

\enddocument
\end